# AGGREGATION FOR GAUSSIAN REGRESSION

By Florentina Bunea,[1] Alexandre B. Tsybakov
and Marten H. Wegkamp[1]

*Florida State University, Université Paris VI and Florida State University*

This paper studies statistical aggregation procedures in the regression setting. A motivating factor is the existence of many different methods of estimation, leading to possibly competing estimators. We consider here three different types of aggregation: model selection (MS) aggregation, convex (C) aggregation and linear (L) aggregation. The objective of (MS) is to select the optimal single estimator from the list; that of (C) is to select the optimal convex combination of the given estimators; and that of (L) is to select the optimal linear combination of the given estimators. We are interested in evaluating the rates of convergence of the excess risks of the estimators obtained by these procedures. Our approach is motivated by recently published minimax results [Nemirovski, A. (2000). Topics in non-parametric statistics. *Lectures on Probability Theory and Statistics (Saint-Flour, 1998)*. *Lecture Notes in Math.* **1738** 85–277. Springer, Berlin; Tsybakov, A. B. (2003). Optimal rates of aggregation. *Learning Theory and Kernel Machines. Lecture Notes in Artificial Intelligence* **2777** 303–313. Springer, Heidelberg]. There exist competing aggregation procedures achieving optimal convergence rates for each of the (MS), (C) and (L) cases separately. Since these procedures are not directly comparable with each other, we suggest an alternative solution. We prove that all three optimal rates, as well as those for the newly introduced (S) aggregation (subset selection), are nearly achieved via a single "universal" aggregation procedure. The procedure consists of mixing the initial estimators with weights obtained by penalized least squares. Two different penalties are considered: one of them is of the BIC type, the second one is a data-dependent $\ell_1$-type penalty.

**1. Introduction.** In this paper we study aggregation procedures and their performance for regression models. Let $\mathcal{D}_n = \{(X_1, Y_1), \ldots, (X_n, Y_n)\}$ be a

Received December 2004; revised September 2006.
[1]Supported in part by NSF Grant DMS-04-06049.
*AMS 2000 subject classifications.* Primary 62G08; secondary 62C20, 62G05, 62G20.
*Key words and phrases.* Aggregation, Lasso estimator, minimax risk, model selection, model averaging, nonparametric regression, oracle inequalities, penalized least squares.

sample of independent random pairs $(X_i, Y_i)$ with

$$(1.1) \qquad Y_i = f(X_i) + W_i, \qquad i = 1, \ldots, n,$$

where $f : \mathcal{X} \to \mathbb{R}$ is an unknown regression function to be estimated, $\mathcal{X}$ is a Borel subset of $\mathbb{R}^d$, the $X_i$'s are fixed elements in $\mathcal{X}$ and the errors $W_i$ are zero mean random variables.

Aggregation of arbitrary estimators in regression models has recently received increasing attention [9, 15, 23, 26, 34, 40, 42, 43, 44, 45]. A motivating factor is the existence of many different methods of estimation, leading to possibly competing estimators $\widehat{f}_1, \ldots, \widehat{f}_M$. A natural idea is then to look for a new, improved estimator $\widetilde{f}$ constructed by combining $\widehat{f}_1, \ldots, \widehat{f}_M$ in a suitable way. Such an estimator $\widetilde{f}$ is called an *aggregate* and its construction is called aggregation.

Three main aggregation problems are model selection (MS) aggregation, convex (C) aggregation and linear (L) aggregation, as first stated by Nemirovski [34]. The objective of (MS) is to select the optimal (in a sense to be defined) single estimator from the list; that of (C) is to select the optimal convex combination of the given estimators; and that of (L) is to select the optimal linear combination of the given estimators.

Aggregation procedures are typically based on sample splitting. The initial sample $\mathcal{D}_n$ is divided into a training sample, used to construct estimators $\widehat{f}_1, \ldots, \widehat{f}_M$, and an independent validation sample, used to learn, that is, to construct $\widetilde{f}$. In this paper we do not consider sample splitting schemes but rather deal with an idealized scheme. We fix the training sample, and thus instead of estimators $\widehat{f}_1, \ldots, \widehat{f}_M$, we have fixed functions $f_1, \ldots, f_M$. A passage to the initial model in our results is straightforward: conditioning on the training sample, we write the inequalities of Theorems 3.1 and 4.1 or, for example, (1.2) below. Then, we take expectations on both sides of these inequalities over the distribution of the whole sample $\mathcal{D}_n$ and interchange the expectation and infimum signs to get bounds containing the risks of the estimators on the right-hand side. The fixed functions $f_1, \ldots, f_M$ can be considered as elements of an (overdetermined) dictionary, see [19], or as "weak learners," see [37], and our results can be interpreted in such contexts as well.

To give precise definitions, denote by $\|g\|_n = \{n^{-1} \sum_{i=1}^{n} g^2(X_i)\}^{1/2}$ the empirical norm of a function $g$ in $\mathbb{R}^d$ and set $\mathsf{f}_\lambda = \sum_{j=1}^{M} \lambda_j f_j$ for any $\lambda = (\lambda_1, \ldots, \lambda_M) \in \mathbb{R}^M$. The performance of an aggregate $\widetilde{f}$ can be judged against the mathematical target

$$(1.2) \qquad \mathbb{E}_f \|\widetilde{f} - f\|_n^2 \leq \inf_{\lambda \in H^M} \|\mathsf{f}_\lambda - f\|_n^2 + \Delta_{n,M},$$

where $\Delta_{n,M} \geq 0$ is a remainder term *independent of $f$* characterizing the price to pay for aggregation, and the set $H^M$ is either the entire $\mathbb{R}^M$ (for linear aggregation), the simplex $\Lambda^M = \{\lambda = (\lambda_1, \ldots, \lambda_M) \in \mathbb{R}^M : \lambda_j \geq 0,$ $\sum_{j=1}^M \lambda_j \leq 1\}$ (for convex aggregation), or the set of all vertices of $\Lambda^M$, except the vertex $(0, \ldots, 0) \in \mathbb{R}^M$ (for model selection aggregation). Here and later $\mathbb{E}_f$ denotes the expectation with respect to the joint distribution of $(X_1, Y_1), \ldots, (X_n, Y_n)$ under model (1.1). The random functions $\mathsf{f}_\lambda$ attaining $\inf_{\lambda \in H^M} \|\mathsf{f}_\lambda - f\|_n^2$ in (1.2) for the three values taken by $H^M$ are called (L), (C) and (MS) oracles, respectively. Note that these minimizers are not estimators since they depend on the true $f$.

We also introduce a fourth type of aggregation, subset selection, or (S) aggregation. For (S) aggregation we fix an integer $D \leq M$ and put $H^M = \Lambda^{M,D}$, where $\Lambda^{M,D}$ denotes the set of all $\lambda \in \mathbb{R}^M$ having at most $D$ nonzero coordinates. Note that (L) aggregation is a special case of subset selection [(S) aggregation] for $D = M$. The literature on subset selection techniques is very large and dates back to [1, 33, 38]. We refer to the recent comprehensive survey [36] for references on methods geared mainly to parametric models. For a review of techniques leading to subset selection in nonparametric settings we refer to [7] and the references therein.

We say that the aggregate $\widetilde{f}$ mimics the (L), (C), (MS) or (S) oracle if it satisfies (1.2) for the corresponding set $H^M$, with the minimal possible price for aggregation $\Delta_{n,M}$. Minimal possible values $\Delta_{n,M}$ for the three problems can be defined via a minimax setting and they are called optimal rates of aggregation [40] and further denoted by $\psi_{n,M}$. Extending the results of [40] obtained in the random design case to the fixed design case, we will show in Sections 3 and 5 that under mild conditions

$$\psi_{n,M} \asymp \begin{cases} M/n, & \text{for (L) aggregation,} \\ \{D \log(1 + M/D)\}/n, & \text{for (S) aggregation,} \\ M/n, & \text{for (C) aggregation, if } M \leq \sqrt{n}, \\ \sqrt{\{\log(1 + M/\sqrt{n})\}/n}, & \text{for (C) aggregation, if } M > \sqrt{n}, \\ (\log M)/n, & \text{for (MS) aggregation.} \end{cases}$$

(1.3)

This implies that linear aggregation has the highest price, (MS) aggregation has the lowest price and convex aggregation occupies an intermediate place. The oracle risks on the right-hand side in (1.2) satisfy a reversed inequality,

$$\inf_{1 \leq j \leq M} \|f_j - f\|_n^2 \geq \inf_{\lambda \in \Lambda^M} \|\mathsf{f}_\lambda - f\|_n^2 \geq \inf_{\lambda \in \mathbb{R}^M} \|\mathsf{f}_\lambda - f\|_n^2,$$

since the sets over which the infima are taken are nested. There is no winner among the three aggregation techniques and the question of how to choose the best among them remains open.

The ideal oracle inequality (1.2) is available only for some special cases. See [13, 15, 27] for (MS) aggregation, [25, 26, 34, 40] for (C) aggregation with $M > \sqrt{n}$ and [40] for (L) aggregation and for (C) aggregation with $M \leq \sqrt{n}$. For more general situations there exist less sharp results of the type

$$(1.4) \qquad \mathbb{E}_f \|\widetilde{f} - f\|_n^2 \leq (1 + \varepsilon) \inf_{\lambda \in H^M} \|\mathrm{f}_\lambda - f\|_n^2 + \Delta_{n,M},$$

where $\varepsilon > 0$ is a constant independent of $f$ and $n$, and $\Delta_{n,M}$ is a remainder term, not necessarily having the same behavior in $n$ and $M$ as the optimal one $\psi_{n,M}$.

Bounds of the type (1.4) in regression problems have been obtained by many authors mainly for the model selection case; see, for example, [4, 5, 7, 8, 9, 10, 11, 12, 15, 23, 28, 30, 32, 42] and the references cited in these works. Most of the papers on model selection treat particular restricted families of estimators, such as orthogonal series estimators, spline estimators, and so forth. There are relatively few results on (MS) aggregation when the estimators are allowed to be arbitrary; see [9, 13, 15, 23, 27, 40, 42, 43, 44, 45]. Various convex aggregation procedures for nonparametric regression have emerged in the last decade. The literature on oracle inequalities of the type (1.2) and (1.4) for the (C) aggregation case is not nearly as large as the one on model selection. We refer to [3, 9, 13, 25, 26, 29, 34, 40, 43, 44, 45]. Finally, linear aggregation procedures are discussed in [13, 34, 40].

Given the existence of competing aggregation procedures achieving either optimal (MS), (C) or (L) bounds, there is an ongoing discussion as to which procedure is the best one. Since this cannot be decided by merely comparing the optimal bounds, we suggest an alternative solution. We show that all three optimal (MS), (C) and (L) bounds can be nearly achieved via a single aggregation procedure. We also show that this procedure leads to near optimal bounds for the newly introduced (S) aggregation, for any subset size $D$. Our answer will thus meet the desiderata of both model (subset) selection and model averaging. The procedures that we suggest for aggregation are based on penalized least squares, with the BIC-type or Lasso ($\ell_1$-type) penalties.

The paper is organized as follows. Section 2 introduces notation and assumptions used throughout the paper. In Section 3 we show that a BIC-type aggregate satisfies inequalities of the form (1.4) with the optimal remainder term $\psi_{n,M}$. We establish the oracle inequalities for all four sets $H^M$ under consideration, hence showing that the BIC-type aggregate achieves simultaneously the (S) [and hence the (L)], the (C) and the (MS) bounds. In Section 4 we study aggregation with the $\ell_1$ penalty and we obtain (1.4) simultaneously for the (S), (C) and (MS) cases, with a remainder term $\Delta_{n,M}$ that differs from the optimal $\psi_{n,M}$ only by a logarithmic factor. We give the corresponding lower bounds for (S), (C) and (MS) aggregation in Section 5, complementing the results obtained for the random design case in [40]. All proofs are deferred to the appendices.

**2. Notation and assumptions.** The following two assumptions on the regression model (1.1) are supposed to be satisfied throughout the paper.

ASSUMPTION (A1). The random variables $W_i$ are independent and Gaussian $N(0, \sigma^2)$.

ASSUMPTION (A2). The functions $f : \mathcal{X} \to \mathbb{R}$ and $f_j : \mathcal{X} \to \mathbb{R}$, $j = 1, \ldots, M$, with $M \geq 2$, belong to the class $\mathcal{F}_0$ of uniformly bounded functions defined by

$$\mathcal{F}_0 \overset{\text{def}}{=} \left\{ g : \mathcal{X} \to \mathbb{R} \,\middle|\, \sup_{x \in \mathcal{X}} |g(x)| \leq L \right\},$$

where $L < \infty$ is a constant that is not necessarily known to the statistician.

The functions $f_j$ can be viewed as estimators of $f$ constructed from a training sample. Here we consider the ideal situation in which they are fixed; we concentrate on learning only. For each $\lambda = (\lambda_1, \ldots, \lambda_M) \in \mathbb{R}^M$, define

$$\mathsf{f}_\lambda(x) = \sum_{j=1}^{M} \lambda_j f_j(x)$$

and let $M(\lambda)$ denote the number of nonzero coordinates of $\lambda$, that is,

$$M(\lambda) = \sum_{j=1}^{M} I_{\{\lambda_j \neq 0\}} = \text{Card } J(\lambda),$$

where $I_{\{\cdot\}}$ denotes the indicator function and $J(\lambda) = \{ j \in \{1, \ldots, M\} : \lambda_j \neq 0 \}$. Furthermore, we introduce the residual sum of squares

$$\widehat{S}(\lambda) = \frac{1}{n} \sum_{i=1}^{n} \{ Y_i - \mathsf{f}_\lambda(X_i) \}^2$$

and the function

$$L(\lambda) = 2 \log \left( \frac{eM}{M(\lambda) \vee 1} \right),$$

for all $\lambda \in \mathbb{R}^M$. The method that we propose is based on aggregating the $f_j$'s via penalized least squares. Given a penalty term $\text{pen}(\lambda)$, the penalized least squares estimator $\widehat{\lambda} = (\widehat{\lambda}_1, \ldots, \widehat{\lambda}_M)$ is defined by

$$(2.1) \qquad \widehat{\lambda} = \underset{\lambda \in \mathbb{R}^M}{\arg \min} \{ \widehat{S}(\lambda) + \text{pen}(\lambda) \},$$

which renders in turn the aggregated estimator

$$(2.2) \qquad \widetilde{f}(x) = \mathsf{f}_{\widehat{\lambda}}(x).$$

Although $\widehat{\lambda}$ is not necessarily unique, all our oracle inequalities hold for any minimizer (2.1). Since the vector $\widehat{\lambda}$ can take any values in $\mathbb{R}^M$, the aggregate $\widetilde{f}$ is not a model selector in the traditional sense, nor is it necessarily a convex combination of the functions $f_j$. Nevertheless, we will show that it mimics the (S), (C) and (MS) oracles when one of the following two penalties is used:

$$(2.3) \qquad \mathrm{pen}(\lambda) = \frac{2\sigma^2}{n}\left\{1 + \frac{2+a}{1+a}\sqrt{L(\lambda)} + \frac{1+a}{a}L(\lambda)\right\}M(\lambda)$$

or

$$(2.4) \qquad \mathrm{pen}(\lambda) = 2\sqrt{2}\sigma\sqrt{\frac{\log M + \log n}{n}}\sum_{j=1}^{M}|\lambda_j|\|f_j\|_n.$$

In (2.3), $a > 0$ is a parameter to be set by the user. The penalty (2.3) can be viewed as a variant of BIC-type penalties [21, 38] since $\mathrm{pen}(\lambda) \sim M(\lambda)$, but the scaling factor here is different and depends on $M(\lambda)$. We also note that in the sequence space model (where the functions $f_1, \ldots, f_M$ are orthonormal with respect to the scalar product induced by the norm $\|\cdot\|_n$), the penalty $\mathrm{pen}(\lambda) \sim M(\lambda)$ leads to $\widehat{\lambda}_j$'s that are hard thresholded values of the $Y_j$'s (see, e.g., [24], page 138). Our penalty (2.3) is not exactly of that form, but it differs from it only in a logarithmic factor.

The penalty (2.4), again in the sequence space model, leads to $\widehat{\lambda}_j$'s that are soft thresholded values of $Y_j$'s. In general models, (2.4) is a weighted $\ell_1$-penalty, with data-dependent weights. Penalized least squares estimators with $\ell_1$-penalty $\mathrm{pen}(\lambda) \sim \sum_{j=1}^{M}|\lambda_j|$ are closely related to basis pursuit [17], to Lasso-type estimators [2, 26, 34, 39] and LARS [20].

Our results show that the BIC-type penalty (2.3) allows optimal aggregation under (A1) and (A2). The $\ell_1$-penalty (2.4) allows near optimal aggregation under somewhat stronger conditions.

## 3. Optimal aggregation with a BIC-type penalty.

In this section we show that the penalized least squares aggregate (2.2) corresponding to the penalty term (2.3) achieves simultaneously the (L), (C) and (MS) bounds of the form (1.4) with the optimal rates $\Delta_{n,M} = \psi_{n,M}$. Consequently, the smallest bound is achieved by our aggregate. The next theorem presents an oracle inequality that implies all three bounds, as well as a bound for (S) aggregation.

THEOREM 3.1. *Assume* (A1) *and* (A2). *Let* $\widetilde{f}$ *be the penalized least squares aggregate defined in* (2.2) *with penalty* (2.3). *Then, for all* $a > 0$ *and all integers* $n \geq 1$ *and* $M \geq 2$,

$$\mathbb{E}_f\|\widetilde{f} - f\|_n^2 \leq (1+a)\inf_{\lambda \in \mathbb{R}^M}\left[\|f_\lambda - f\|_n^2 + \frac{\sigma^2}{n}\left\{5 + \frac{2+3a}{a}L(\lambda)\right\}M(\lambda)\right]$$

(3.1)
$$+ \frac{\sigma^2}{n} \frac{6(1+a)^2}{a(e-1)}.$$

The proof is given in Appendix A.

COROLLARY 3.2. *Under the conditions of Theorem* 3.1, *there exists a constant $C > 0$ such that for all $a > 0$ and all integers $n \geq 1$ and $M \geq 2$ and $D \leq M$, the following upper bounds for $R_{M,n} \stackrel{\text{def}}{=} \mathbb{E}_f \|\tilde{f} - f\|_n^2$ hold:*

(3.2)  $$R_{M,n} \leq (1+a) \inf_{1 \leq j \leq M} \|f_j - f\|_n^2 + C(1 + a + a^{-1})\sigma^2 \frac{\log M}{n},$$

(3.3)  $$R_{M,n} \leq (1+a) \inf_{\lambda \in \Lambda^{M,D}} \|\mathbf{f}_\lambda - f\|_n^2 + C(1 + a + a^{-1})\sigma^2 \frac{D}{n} \log\left(\frac{M}{D} + 1\right),$$

(3.4)  $$R_{M,n} \leq (1+a) \inf_{\lambda \in \mathbb{R}^M} \|\mathbf{f}_\lambda - f\|_n^2 + C(1 + a + a^{-1})\sigma^2 \frac{M}{n},$$

(3.5)  $$R_{M,n} \leq (1+a) \inf_{\lambda \in \Lambda^M} \|\mathbf{f}_\lambda - f\|_n^2 + C(1 + a + a^{-1})(L^2 + \sigma^2)\psi_n^C(M),$$

*where*

$$\psi_n^C(M) = \begin{cases} M/n, & \text{if } M \leq \sqrt{n}, \\ \sqrt{\{\log(eM/\sqrt{n})\}/n}, & \text{if } M > \sqrt{n}. \end{cases}$$

The proof is given in Appendix A.

Note that, along with the bounds of Corollary 3.2, Theorem 3.1 implies a trivial (constant) bound on $R_{M,n}$. In fact, the infimum over $\lambda \in \mathbb{R}^M$ in (3.1) is smaller than the value of the expression in square brackets for $\lambda = 0$, which together with Assumption (A2) yields $R_{M,n} \leq (1+a)L^2 + \frac{\sigma^2}{n}\frac{6(1+a)^2}{a(e-1)}$. Therefore, the remainder terms in (3.2)–(3.5) can be replaced by their truncated versions (truncation by a constant).

A variant of Theorem 3.1 for regression with random design $X_1, \ldots, X_n$ can be found in [14].

REMARK 1. The variance $\sigma^2$ is usually unknown and we need to substitute an estimate in the penalty (2.3). We consider the situation described in the introduction where the functions $f_j$ are estimators based on an independent (training) data set $\mathcal{D}'_\ell$ that consists of observations $(X'_j, Y'_j)$ following (1.1). Let $\hat{\sigma}^2$ be an estimate of $\sigma^2$ based on $\mathcal{D}'_\ell$ only. We write $\mathbb{E}_f^{(1)}$ ($\mathbb{E}_f^{(2)}$) for expectation with respect to $\mathcal{D}'_\ell$ ($\mathcal{D}_n$), and let $\mathbb{E}_f = \mathbb{E}_f^{(1)}\mathbb{E}_f^{(2)}$ be the product expectation. Let $\hat{f}$ be the aggregate corresponding to the penalty (2.3) with

$\sigma^2$ replaced by $\widehat{\sigma}^2$. Note that

$$\mathbb{E}_f \|\widehat{f} - f\|_n^2 = \mathbb{E}_f^{(1)} \mathbb{E}_f^{(2)} \|\widehat{f} - f\|_n^2 I_{\{2\widehat{\sigma}^2 \geq \sigma^2\}} + \mathbb{E}_f^{(1)} \mathbb{E}_f^{(2)} \|\widehat{f} - f\|_n^2 I_{\{2\widehat{\sigma}^2 < \sigma^2\}}.$$

Inspection of the proof of Theorem 3.1 shows that $\mathbb{E}_f^{(2)} \|\widehat{f} - f\|_n^2 I_{\{2\widehat{\sigma}^2 \geq \sigma^2\}}$ can be bounded simply by the right-hand side of (3.1) with $\sigma^2$ substituted by $2\widehat{\sigma}^2$, as Theorem 3.1 holds for any penalty term larger than (2.3). Consequently we find

$$\mathbb{E}_f \|\widehat{f} - f\|_n^2 I_{\{2\widehat{\sigma}^2 \geq \sigma^2\}}$$

$$\leq \frac{2\mathbb{E}_f^{(1)}\widehat{\sigma}^2}{n} \frac{6(1+a)^2}{a(e-1)}$$

$$+ (1+a) \inf_{\lambda \in \mathbb{R}^M} \left[ \|\mathsf{f}_\lambda - f\|_n^2 + \frac{2\mathbb{E}_f^{(1)}\widehat{\sigma}^2}{n} \left\{ 5 + \frac{2+3a}{a} L(\lambda) \right\} M(\lambda) \right].$$

Next, we observe that $\mathbb{E}_f^{(2)} \|\widehat{f} - f\|_n^2 \leq 6\sigma^2 + 2L^2$. For this, we use the reasoning leading to (A.5) in the proof of Theorem 4.1, in which we replace $I_{A^c}$ by 1 throughout. Notice that this argument holds for any positive penalty term $\mathrm{pen}(\lambda)$ such that $\mathrm{pen}(\lambda_0) = 0$ with $\lambda_0 = (0, \ldots, 0)$, and hence it holds for the penalty term used here. Thus

$$\mathbb{E}_f \|\widehat{f} - f\|_n^2 I_{\{2\widehat{\sigma}^2 < \sigma^2\}} \leq (6\sigma^2 + 2L^2) \mathbb{P}_f^{(1)} \{2\widehat{\sigma}^2 < \sigma^2\}.$$

Combining the three displays above we see that $\widehat{f}$ achieves a bound similar to (3.1) if the estimator $\widehat{\sigma}^2$ satisfies $\mathbb{P}_f^{(1)}\{2\widehat{\sigma}^2 < \sigma^2\} \leq c_1/n$ and $\mathbb{E}_f^{(1)}\widehat{\sigma}^2 \leq c_2\sigma^2$ for some finite constants $c_1, c_2$. Since the sample variance of $Y_i'$ from the training sample $\mathcal{D}_\ell'$, with $\ell \geq cn$ for some positive constant $c$, meets both requirements, it can always play the role of $\widehat{\sigma}^2$.

**4. Near optimal aggregation with a data-dependent $\ell_1$ penalty.** In this section we show that the penalized least squares aggregate (2.2) using a penalty of the form (2.4) achieves simultaneously the (MS), (C), (L) and (S) bounds of the form (1.4) with near optimal rates $\Delta_{n,M}$. We will use the following additional assumption.

ASSUMPTION (A3). Define the matrices

$$\Psi_n = \left( \frac{1}{n} \sum_{i=1}^n f_j(X_i) f_{j'}(X_i) \right)_{1 \leq j,j' \leq M},$$

$$\mathrm{diag}(\Psi_n) = \mathrm{diag}(\|f_1\|_n^2, \ldots, \|f_M\|_n^2).$$

There exists $\kappa = \kappa_{n,M} > 0$ such that the matrix $\Psi_n - \kappa \, \mathrm{diag}(\Psi_n)$ is positive semidefinite for any given $n \geq 1$, $M \geq 2$.

The next theorem presents an oracle inequality similar to the one of Theorem 3.1.

THEOREM 4.1. *Assume* (A1), (A2) *and* (A3). *Let* $\widetilde{f}$ *be the penalized least squares aggregate defined by* (2.2) *with penalty* (2.4). *Then, for all* $\varepsilon > 0$, *and all integers* $n \geq 1$, $M \geq 2$, *we have*

$$
\mathbb{E}_f \|\widetilde{f} - f\|_n^2
$$

$$
\begin{aligned}
(4.1) \qquad &\leq \inf_{\lambda \in \mathbb{R}^M} \left\{ (1+\varepsilon)\|\mathsf{f}_\lambda - f\|_n^2 + 8\left(4 + \varepsilon + \frac{4}{\varepsilon}\right)\frac{\sigma^2}{\kappa}\frac{\log M + \log n}{n}M(\lambda) \right\} \\
&\quad + \frac{4L^2 + 12\sigma^2}{n\sqrt{\pi(\log M + \log n)}} + 6\sigma^2\sqrt{\frac{n+2}{n}}\exp\left(-\frac{n}{16}\right).
\end{aligned}
$$

The proof is given in Appendix A.

COROLLARY 4.2. *Let the assumptions of Theorem* 4.1 *be satisfied. Then there exists a constant* $C = C(L^2, \sigma^2, \kappa) > 0$ *such that for all* $\varepsilon > 0$ *and for all integers* $n \geq 1$, $M \geq 2$ *and* $D \leq M$,

$$
\mathbb{E}_f\|\widetilde{f} - f\|_n^2 \leq (1+\varepsilon)\inf_{1 \leq j \leq M}\|f_j - f\|_n^2 + C(1+\varepsilon+\varepsilon^{-1})\frac{\log(M \vee n)}{n},
$$

$$
\mathbb{E}_f\|\widetilde{f} - f\|_n^2 \leq (1+\varepsilon)\inf_{\lambda \in \Lambda^{M,D}}\|\mathsf{f}_\lambda - f\|_n^2 + C(1+\varepsilon+\varepsilon^{-1})\frac{D\log(M \vee n)}{n},
$$

$$
\mathbb{E}_f\|\widetilde{f} - f\|_n^2 \leq (1+\varepsilon)\inf_{\lambda \in \mathbb{R}^M}\|\mathsf{f}_\lambda - f\|_n^2 + C(1+\varepsilon+\varepsilon^{-1})\frac{M\log(M \vee n)}{n},
$$

$$
\mathbb{E}_f\|\widetilde{f} - f\|_n^2 \leq (1+\varepsilon)\inf_{\lambda \in \Lambda^M}\|\mathsf{f}_\lambda - f\|_n^2 + C(1+\varepsilon+\varepsilon^{-1})\overline{\psi_n^C}(M),
$$

*where*

$$
\overline{\psi_n^C}(M) = \begin{cases} (M\log n)/n, & \text{if } M \leq \sqrt{n}, \\ \sqrt{(\log M)/n}, & \text{if } M > \sqrt{n}. \end{cases}
$$

PROOF. The argument is similar to that of the proof of Corollary 3.2. □

REMARK 2. Using the same reasoning as in Remark 1, we can replace $\sigma^2$ in the penalty term by twice the sample variance of $Y_i'$ from the training sample $\mathcal{D}_\ell'$.

REMARK 3. Inspection of the proofs shows that the constants $C = C(L^2, \sigma^2, \kappa)$ in Corollary 4.2 have the form $C = A_1 + A_2/\kappa$, where $A_1$ and $A_2$ are constants independent of $\kappa$. In general, $\kappa$ may depend on $n$ and $M$.

However, if $\kappa > c$ for some constant $c > 0$, independent of $n$ and $M$, as discussed in Remarks 3 and 4 below, the rates of aggregation given in Corollary 4.2 are near optimal, up to logarithmic factors. They are exactly optimal [cf. (1.3) and the lower bounds of the next section] for some configurations of $n, M$: for (MS)-aggregation if $n^{a'} \le M \le n^a$, and for (C)-aggregation if $n^{1/2} \le M \le n^a$, where $0 < a' < a < \infty$.

REMARK 4.   If $\xi_{\min}$, the smallest eigenvalue of the matrix $\Psi_n$, is positive, Assumption (A3) is satisfied for $\kappa = \xi_{\min}/L^2$. In a standard parametric regression context where $M$ is fixed and $\Psi_n/n$ converges to a nonsingular $M \times M$ matrix, we have that $\xi_{\min} > c$ (and therefore $\kappa > c/L^2$) for $n$ large enough and for some $c > 0$ independent of $M$ and $n$.

REMARK 5.   Assumption (A3) is trivially satisfied with $\kappa = 1$ if $\Psi_n$ is a diagonal matrix (note that zero diagonal entries are not excluded). An example illustrating this situation is related to the orthogonal series nonparametric regression: $M = M_n$ is allowed to converge to $\infty$ as $n \to \infty$ and the basis functions $f_j$ are orthogonal with respect to the empirical norm. Another example is related to sequence space models, where the $f_j$ are estimators constructed from nonintersecting blocks of coefficients. Aggregation of such mutually orthogonal estimators can be used to achieve adaptation (cf., e.g., [34]). Note that Assumption (A3) does not exclude the matrices $\Psi_n$ whose ordered eigenvalues can be arbitrarily close to 0 as $M \to \infty$. The last property is characteristic for sequence space representation of statistical inverse problems: there $\Psi_n$ is diagonal, with $M = M_n \to \infty$, as $n \to \infty$, and with the $j$th eigenvalue tending to 0 as $j \to \infty$. For such matrices $\Psi_n$ Assumption (A3) holds with $\kappa = 1$, so that the oracle inequality of Theorem 4.1 is invariant with respect to the speed of decrease of the eigenvalues.

REMARK 6.   The bounds of Corollary 4.2 can be written with the remainder terms truncated at a constant level (cf. an analogous remark after Corollary 3.2). Thus, for $M > n$ the (L) bounds become trivial.

However, for $M > n$ an oracle bound of the type (4.1) is still meaningful if $f$ is sparse, that is, can be well approximated by relatively few (less than $n$) functions $f_j$. This is illustrated by the next theorem where Assumption (A3) is replaced by a local mutual coherence property of the matrix $\Psi_n$, relaxing the mutual coherence condition suggested in [19]. Let

$$\rho_n(i, j) = \frac{\langle f_i, f_j \rangle_n}{\|f_i\|_n \|f_j\|_n}$$

denote the "correlation" between two elements $f_i$ and $f_j$. We will assume that the values $\rho_M(i, j)$ with $i \ne j$ are relatively small, for $i \in J(\lambda)$, $\lambda \in \mathbb{R}^M$.

Set

$$\rho(\lambda) = \max_{i \in J(\lambda)} \max_{j > i} |\rho_n(i,j)|.$$

THEOREM 4.3. *Assume* (A1) *and* (A2)*. Let* $\widetilde{f}$ *be the penalized least squares aggregate defined by* (2.2) *with penalty*

$$\mathrm{pen}(\lambda) = 4\sqrt{2}\sigma\sqrt{\frac{\log M + \log n}{n}} \sum_{j=1}^{M} |\lambda_j| \|f_j\|_n.$$

*Then, for all* $\varepsilon > 0$ *and all integers* $n \geq 1$, $M \geq 2$, *we have*

$$\mathbb{E}_f \|\widetilde{f} - f\|_n^2 \leq \inf \left\{ (1+\varepsilon)\|\mathsf{f}_\lambda - f\|_n^2 + 32\left(4 + \varepsilon + \frac{4}{\varepsilon}\right)\sigma^2 \frac{\log M + \log n}{n} M(\lambda) \right\}$$
$$+ \frac{4L^2 + 12\sigma^2}{n\sqrt{\pi(\log M + \log n)}} + 6\sigma^2 \sqrt{\frac{n+2}{n}} \exp\left(-\frac{n}{16}\right),$$

*where the infimum is taken over all* $\lambda \in \mathbb{R}^M$ *such that* $32\rho(\lambda)M(\lambda) \leq 1$.

The proof is given in Appendix A.

In particular, if $f$ has a sparse representation $f = \mathsf{f}_{\lambda^*}$ for some $\lambda^* \in \mathbb{R}^M$ with $32\rho(\lambda^*)M(\lambda^*) \leq 1$, there exists a constant $C = C(L^2, \sigma^2) < \infty$ such that

$$\mathbb{E}_f \|\widetilde{f} - f\|_n^2 \leq C(\log M + \log n)\frac{M(\lambda^*)}{n}$$

for all $n \geq 1$ and $M \geq 2$. Even for $M > n$, this bound is meaningful if $M(\lambda^*) \ll n$.

Note that in Theorem 4.3 the correlations $\rho_n(i,j)$ with $i, j \notin J(\lambda)$ can take arbitrary values in $[-1, 1]$. Such $\rho_n(i,j)$ constitute the overwhelming majority of the elements of the correlation matrix if $J(\lambda)$ is a set of small cardinality, $M(\lambda) \ll M$.

REMARK 7. An attractive feature of the $\ell_1$-penalized aggregation is its computational feasibility. Clearly, the criterion in (2.1) with penalties as in Theorems 4.1 and 4.3 is convex in $\lambda$. One can therefore use techniques of convex optimization to compute the aggregates. We refer, for instance, to [20, 35] for detailed analysis of such optimization problems and fast algorithms.

REMARK 8. We refer to Theorem 2.1 in [35] for conditions under which the penalized least squares aggregate is unique. Typically, for $M > n$ the solution is not unique, but a convex combination of solutions is itself a solution. Our results hold for any element of such a convex set of solutions.

**5. Lower bounds.** In this section we provide lower bounds showing that the remainder terms in the upper bounds obtained in the previous sections are optimal or near optimal. For regression with random design and the $L_2(\mathbb{R}^d, \mu)$-risks, such lower bounds for aggregation with optimal rates $\psi_{n,M}$ as given in (1.3) were established in [40]. The next theorem extends them to aggregation for the regression model with fixed design. Furthermore, we state these bounds in a more general form, considering not only the expected squared risks, but also other loss functions, and instead of the (L) aggregation lower bound, we provide the more general (S) aggregation bound.

Let $w:\mathbb{R} \to [0,\infty)$ be a *loss function*, that is, a monotone nondecreasing function satisfying $w(0) = 0$ and $w \not\equiv 0$.

THEOREM 5.1. *Let the integers $n, M, D$ be such that $2 \leq M \leq n$, and let $X_1,\ldots,X_n$ be distinct points. Assume that $H^M$ is either the simplex $\Lambda^M$ [for the (C) aggregation case], the set of vertices of $\Lambda^M$, except the vertex $(0,\ldots,0) \in \mathbb{R}^M$ [for the (MS) aggregation case], or the set $\Lambda^{M,D}$ [for the (S) aggregation case]. Let the corresponding $\psi_{n,M}$ be given by (1.3) and, for (S) aggregation, assume that $M \log(M/D+1) \leq n$ and $M \geq D$. Then there exist $f_1,\ldots,f_M \in \mathcal{F}_0$ such that*

$$(5.1) \qquad \inf_{T_n} \sup_{f \in \mathcal{F}_0} \mathbb{E}_f w\left[\psi_{n,M}^{-1}\left(\|T_n - f\|_n^2 - \inf_{\lambda \in H^M}\|f_\lambda - f\|_n^2\right)\right] \geq c,$$

*where $\inf_{T_n}$ denotes the infimum over all estimators and the constant $c > 0$ does not depend on $n, M$ and $D$.*

The proof is given in Appendix A.

Setting $w(u) = u$ in Theorem 5.1, we get the lower bounds for expected squared risks showing optimality or near optimality of the remainder terms in the oracle inequalities of Corollaries 3.2 and 4.2. The choice of $w(u) = I\{u > a\}$ with some fixed $a > 0$ leads to the lower bounds for probabilities showing near optimality of the remainder terms in the corresponding upper bounds "in probability" obtained in [14].

## APPENDIX A: PROOFS

**A.1. Proof of Theorem 3.1.** Let $\Lambda_m$ be the set of $\lambda \in \mathbb{R}^M$ with exactly $m$ nonzero coefficients, $\Lambda_m = \{\lambda \in \mathbb{R}^M : M(\lambda) = m\}$. Let $J_{m,k}$, $k = 1,\ldots,\binom{M}{m}$, be all the subsets of $\{1,\ldots,M\}$ of cardinality $m$ and define

$$\Lambda_{m,k}(\lambda) = \{\lambda = (\lambda_1,\ldots,\lambda_M) \in \Lambda_m : \lambda_j \neq 0 \Leftrightarrow j \in J_{m,k}\}.$$

The collection $\{\Lambda_{m,k} : 1 \leq k \leq \binom{M}{m}\}$ forms a partition of the set $\Lambda_m$. Observe that

$$\inf_{\lambda \in \mathbb{R}^M}\{\widehat{S}(\lambda) + \text{pen}(\lambda)\} = \inf_{0 \leq m \leq M} \inf_{1 \leq k \leq \binom{M}{m}} \inf_{\lambda \in \Lambda_{m,k}} \{\widehat{S}(\lambda) + \text{pen}(\lambda)\}.$$

Here the penalty pen($\lambda$) is defined in (2.3), and it takes a constant value on each of the sets $\Lambda_m$ as $M(\lambda) = m$ and $L(\lambda) = L_m \equiv 2\ln(eM/(m \vee 1))$ for all $\lambda \in \Lambda_m$. We now apply [11], Theorem 2, choosing there the parameters $\theta = a/(1+a)$ and $K = 2$. This yields

$$\mathbb{E}_f \|\tilde{f} - f\|_n^2 \le (1+a) \inf_{0 \le m \le M} \inf_{1 \le k \le \binom{M}{m}} \left\{ \inf_{\lambda \in \Lambda_{m,k}} \|f_\lambda - f\|_n^2 + \text{pen}(\lambda) - \frac{m\sigma^2}{n} \right\}$$

$$+ \frac{(1+a)^2}{a} \frac{\sigma^2}{n} \Sigma \left\{ \frac{(2+a)^2}{(1+a)^2} + 2 \right\},$$

where $\Sigma = \sum_{m=1}^{M} \binom{M}{m} \exp(-mL_m)$. Using the crude bound $\binom{M}{m} \le (eM/m)^m$ (see, e.g., [18], page 218), we get

$$\Sigma \le \sum_{m=1}^{M} \left( \frac{eM}{m} \right)^{-m} \le \sum_{m=1}^{M} e^{-m} \le \frac{1}{e-1}.$$

For all $\lambda \in \Lambda_m$, we have

$$n\,\text{pen}(\lambda) - m\sigma^2 = \sigma^2 m \left( 1 + 2\frac{2+a}{1+a}\sqrt{L_m} + 2\frac{1+a}{a}L_m \right)$$

$$\le \sigma^2 m \left( 5 + \frac{2+3a}{a}L_m \right).$$

Consequently we find

$$\mathbb{E}_f \|\tilde{f} - f\|_n^2 \le (1+a)$$

$$\times \inf_{0 \le m \le M} \inf_{1 \le k \le \binom{M}{m}} \left\{ \inf_{\lambda \in \Lambda_{m,k}} \|f - f_\lambda\|_n^2 + \frac{\sigma^2 m}{n} \left( 5 + \frac{2+3a}{a}L_m \right) \right\}$$

$$+ \frac{6(1+a)^2}{a(e-1)} \frac{\sigma^2}{n}$$

$$= (1+a) \inf_{\lambda \in \mathbb{R}^M} \left[ \|f - f_\lambda\|_n^2 + \frac{\sigma^2 M(\lambda)}{n} \left\{ 5 + \frac{2+3a}{a}L(\lambda) \right\} \right]$$

$$+ \frac{6(1+a)^2}{a(e-1)} \frac{\sigma^2}{n},$$

which proves the result.

### A.2. Proof of Corollary 3.2.

PROOF OF (3.2). Since the infimum on the right-hand side of (3.1) is taken over all $\lambda \in \mathbb{R}^M$, the bound easily follows by restricting the minimization to the set of the $M$ vertices $(1, 0, \ldots, 0)$, $(0, 1, 0, \ldots, 0), \ldots, (0, \ldots, 0, 1)$ of $\Lambda^M$. $\quad\square$

PROOF OF (3.3) AND (3.4). The (S) bound (3.3) easily follows from (3.1) by restricting the minimization to $\Lambda^{M,D}$. In fact, for $\lambda \in \Lambda^{M,D}$ we have $M(\lambda) = D$ and $L(\lambda) = 2\log(eM/D) \le 6\log(M/D+1)$. The (L) bound (3.4) is a special case of (3.3) for $D = M$. □

PROOF OF (3.5). For $M \le \sqrt{n}$ the result follows from (3.4). Assume now that $M > \sqrt{n}$ and let $m \ge 1$ be the smallest integer greater than or equal to

$$x_{n,M} = \sqrt{n}/(2\sqrt{\log(eM/\sqrt{n})}).$$

Clearly, $x_{n,M} \le m \le x_{n,M} + 1 \le M$. First, consider the case $m \ge 1$. Denote by $\mathcal{C}$ the set of functions $h$ of the form

$$h(x) = \frac{1}{m}\sum_{j=1}^{M} k_j f_j(x), \qquad k_j \in \{0, 1, \ldots, m\}, \ \sum_{j=1}^{M} k_j \le m.$$

The following approximation result can be obtained by the "Maurey argument" (see, e.g., [6], Lemma 1 or [34], pages 192 and 193):

$$(A.1) \qquad \min_{g \in \mathcal{C}} \|g - f\|_n^2 \le \min_{\lambda \in \Lambda^M} \|\mathsf{f}_\lambda - f\|_n^2 + \frac{L^2}{m}.$$

For completeness, we give the proof of (A.1) in Appendix B. Since $M(\lambda) \le m \le x_{n,M} + 1$ for the vectors $\lambda$ corresponding to $g \in \mathcal{C}$, and since $x \mapsto x\log(\frac{eM}{x})$ is increasing for $1 \le x \le M$, we get from (3.1)

$$\mathbb{E}_f \|\widetilde{f} - f\|_n^2 \le \inf_{g \in \mathcal{C}} \left\{ C_1 \|g - f\|_n^2 + C_2 \frac{x_{n,M}+1}{n}\log\left(\frac{eM}{x_{n,M}+1}\right) \right\} + \frac{C_3}{n}$$

with $C_1 = 1+a, C_2 = C_2'(1+a+1/a)\sigma^2$ and $C_3 = C_3'(1+a+1/a)\sigma^2$, where $C_2' > 0$ and $C_3' > 0$ are absolute constants. Using this inequality, (A.1) and the fact that $m \ge x_{n,M}$, we obtain

$$\mathbb{E}_f \|\widetilde{f} - f\|_n^2 \le C_1 \inf_{\lambda \in \Lambda^M} \|\mathsf{f}_\lambda - f\|_n^2 + C_1 \frac{L^2}{x_{n,M}} + C_2 \frac{x_{n,M}+1}{n}\log\left(\frac{eM}{x_{n,M}}\right) + \frac{C_3}{n}.$$

Since, clearly, $n^{-1} \le \psi_n^C(M)$, to complete the proof of (3.5) it remains to note that

$$\log\left(\frac{eM}{x_{n,M}}\right) = \log\left(\frac{2eM}{\sqrt{n}}\sqrt{\log\left(\frac{eM}{\sqrt{n}}\right)}\right) \le 3\log\left(\frac{eM}{\sqrt{n}}\right),$$

in view of the elementary inequality $\log(2y\sqrt{\log(y)}) \le 3\log(y)$, for all $y \ge e$. □

**A.3. Proof of Theorem 4.1.** We begin as in [31]. First we define

$$r_n = 2\sqrt{2}\sigma\sqrt{\frac{\log M + \log n}{n}}$$

and $r_{n,j} = r_n\|f_j\|_n$. By definition, $\widetilde{f} = \mathsf{f}_{\widehat{\lambda}}$ satisfies

$$\widehat{S}(\widehat{\lambda}) + \sum_{j=1}^{M} r_{n,j}|\widehat{\lambda}_j| \leq \widehat{S}(\lambda) + \sum_{j=1}^{M} r_{n,j}|\lambda_j|$$

for all $\lambda \in \mathbb{R}^M$, which we may rewrite as

$$\|\widetilde{f} - f\|_n^2 + \sum_{j=1}^{M} r_{n,j}|\widehat{\lambda}_j| \leq \|\mathsf{f}_\lambda - f\|_n^2 + \sum_{j=1}^{M} r_{n,j}|\lambda_j| + \frac{2}{n}\sum_{i=1}^{n} W_i(\widetilde{f} - \mathsf{f}_\lambda)(X_i).$$

We define the random variables $V_j = \frac{1}{n}\sum_{i=1}^{n} f_j(X_i)W_i$, $1 \leq j \leq M$, and the event

$$A = \bigcap_{j=1}^{M}\{2|V_j| \leq r_{n,j}\}.$$

The normality Assumption (A1) on $W_i$ implies that $\sqrt{n}V_j \sim N(0, \sigma^2\|f_j\|_n^2)$, $1 \leq j \leq M$. Applying the union bound followed by the standard tail bound for the $N(0,1)$ distribution, we find

$$(\text{A.2}) \qquad \mathbb{P}(A^c) \leq \sum_{j=1}^{M} \mathbb{P}\{\sqrt{n}|V_j| > \sqrt{n}r_{n,j}/2\} \leq \frac{1}{n\sqrt{\pi(\log M + \log n)}}.$$

Then, on the set $A$, we find

$$\frac{2}{n}\sum_{i=1}^{n} W_i(\widetilde{f} - \mathsf{f}_\lambda)(X_i) = 2\sum_{j=1}^{M} V_j(\widehat{\lambda}_j - \lambda_j) \leq \sum_{j=1}^{M} r_{n,j}|\widehat{\lambda}_j - \lambda_j|$$

and therefore, still on the set $A$,

$$\|\widetilde{f} - f\|_n^2 \leq \|\mathsf{f}_\lambda - f\|_n^2 + \sum_{j=1}^{M} r_{n,j}|\widehat{\lambda}_j - \lambda_j| + \sum_{j=1}^{M} r_{n,j}|\lambda_j| - \sum_{j=1}^{M} r_{n,j}|\widehat{\lambda}_j|.$$

Recall that $J(\lambda)$ denotes the set of indices of the nonzero elements of $\lambda$, and $M(\lambda) = \text{Card } J(\lambda)$. Rewriting the right-hand side of the previous display, we find, on the set $A$,

$$\|\widetilde{f} - f\|_n^2 \leq \|\mathsf{f}_\lambda - f\|_n^2 + \left(\sum_{j=1}^{M} r_{n,j}|\widehat{\lambda}_j - \lambda_j| - \sum_{j\notin J(\lambda)} r_{n,j}|\widehat{\lambda}_j|\right)$$

$$\text{(A.3)} \qquad + \left( -\sum_{j \in J(\lambda)} r_{n,j} |\widehat{\lambda}_j| + \sum_{j \in J(\lambda)} r_{n,j} |\lambda_j| \right)$$

$$\leq \|\mathsf{f}_\lambda - f\|_n^2 + 2 \sum_{j \in J(\lambda)} r_{n,j} |\widehat{\lambda}_j - \lambda_j|$$

by the triangle inequality and the fact that $\lambda_j = 0$ for $j \notin J(\lambda)$. By Assumption (A3), we have

$$\sum_{j \in J(\lambda)} r_{n,j}^2 |\widehat{\lambda}_j - \lambda_j|^2 \leq r_n^2 \sum_{j=1}^M \|f_j\|_n^2 |\widehat{\lambda}_j - \lambda_j|^2 = r_n^2 (\widehat{\lambda} - \lambda)' \operatorname{diag}(\Psi_n)(\widehat{\lambda} - \lambda)$$

$$\leq r_n^2 \kappa^{-1} (\widehat{\lambda} - \lambda)' \Psi_n (\widehat{\lambda} - \lambda) = r_n^2 \kappa^{-1} \|\widetilde{f} - \mathsf{f}_\lambda\|_n^2.$$

Combining this with the Cauchy–Schwarz and triangle inequalities, we find further that, on the set $A$,

$$\|\widetilde{f} - f\|_n^2 \leq \|\mathsf{f}_\lambda - f\|_n^2 + 2 \sum_{j \in J(\lambda)} r_{n,j} |\widehat{\lambda}_j - \lambda_j|$$

$$\text{(A.4)} \qquad \leq \|\mathsf{f}_\lambda - f\|_n^2 + 2 r_n \sqrt{M(\lambda)/\kappa} (\|\widetilde{f} - f\|_n + \|\mathsf{f}_\lambda - f\|_n).$$

Inequality (A.4) is of the simple form $v^2 \leq c^2 + vb + cb$ with $v = \|\widetilde{f} - f\|_n$, $b = 2 r_n \sqrt{M(\lambda)/\kappa}$ and $c = \|\mathsf{f}_\lambda - f\|_n$. After applying the inequality $2xy \leq x^2/\alpha + \alpha y^2$ $(x, y \in \mathbb{R}, \alpha > 0)$ twice, to $2bc$ and $2bv$, we easily find $v^2 \leq v^2/(2\alpha) + \alpha b^2 + (2\alpha+1)/(2\alpha)c^2$, whence $v^2 \leq a/(a-1)\{b^2(a/2) + c^2(a+1)/a\}$ for $a = 2\alpha > 1$. Recalling that (A.4) is valid on the set $A$, we now get that

$$\mathbb{E}_f[\|\widetilde{f} - f\|_n^2 I_A] \leq \inf_{\lambda \in \mathbb{R}^M} \left\{ \frac{a+1}{a-1} \|\mathsf{f}_\lambda - f\|_n^2 + \frac{2a^2}{\kappa(a-1)} r_n^2 M(\lambda) \right\} \qquad \forall\, a > 1.$$

It remains to bound $\mathbb{E}_f \|\widetilde{f} - f\|_n^2 I_{A^c}$. Writing $\|W\|_n^2 = n^{-1} \sum_{i=1}^n W_i^2$ and using the inequality $(x+y)^2 \leq 2x^2 + 2y^2$, we find that

$$\mathbb{E}_f \|\widetilde{f} - f\|_n^2 I_{A^c} \leq 2 \mathbb{E}_f \widehat{S}(\widetilde{f}) I_{A^c} + 2 \mathbb{E}_f \|W\|_n^2 I_{A^c}.$$

Next, since $\operatorname{pen}(\widetilde{\lambda}) \geq 0$ and by the definition of $\widetilde{f}$, for $\lambda_0 = (0, \ldots, 0)' \in \mathbb{R}^M$,

$$\mathbb{E}_f \widehat{S}(\widetilde{f}) I_{A^c} \leq \mathbb{E}_f \{\widehat{S}(\widetilde{f}) + \operatorname{pen}(\widetilde{\lambda})\} I_{A^c} \leq \mathbb{E}_f \{\widehat{S}(f_{\lambda_0}) + \operatorname{pen}(\lambda_0)\} I_{A^c}$$

$$= \mathbb{E}_f \widehat{S}(f_{\lambda_0}) I_{A^c} \leq 2 L^2 \mathbb{P}(A^c) + 2 \mathbb{E}_f \|W\|_n^2 I_{A^c},$$

whence

$$\text{(A.5)} \qquad \mathbb{E}_f \|\widetilde{f} - f\|_n^2 I_{A^c} \leq 4 L^2 \mathbb{P}(A^c) + 6 \mathbb{E}_f \|W\|_n^2 I_{A^c}.$$

In order to bound the last term on the right-hand side, we introduce the event $B = \{\frac{1}{n}\sum_{i=1}^{n} W_i^2 \leq 2\sigma^2\}$. Using Lemma B.2 from Appendix B with $d = n$, we get

$$\mathbb{P}\{B^c\} = \mathbb{P}\{Z_n^2 - n > \sqrt{2n}\sqrt{n/2}\} \leq \exp(-n/8).$$

Observe further that $\mathbb{E}_f \|W\|_n^2 I_{A^c} \leq 2\sigma^2 \mathbb{P}\{A^c\} + \mathbb{E}_f \|W\|_n^2 I_{B^c}$ and by the Cauchy–Schwarz inequality we find

$$\mathbb{E}_f \|W\|_n^2 I_{B^c} \leq (\mathbb{E}_f \|W\|_n^4)^{1/2} \exp(-n/16) = \left(\frac{3\sigma^4}{n} + \frac{n-1}{n}\sigma^4\right)^{1/2} \exp(-n/16).$$

Collecting all these bounds, and using the bound (A.2) on $\mathbb{P}\{A^c\}$, we obtain

$$\mathbb{E}_f \|\widetilde{f} - f\|_n^2 I_{A^c} \leq 4L^2 \mathbb{P}(A^c) + 6\mathbb{E}_f \|W\|_n^2 I_{A^c}$$

$$\leq \frac{4L^2 + 12\sigma^2}{n\sqrt{\pi(\log M + \log n)}} + 6\sigma^2 \sqrt{\frac{n+2}{n}} \exp(-n/16).$$

The proof of the theorem is complete by taking $\varepsilon = 2/(a-1)$.

**A.4. Proof of Theorem 4.3.**    First, notice that by definition of $\widetilde{f}$ and of the penalty $\mathrm{pen}(\lambda) = 2\sum_{j=1}^{M} r_{n,j}|\lambda_j|$,

$$\|\widetilde{f} - f\|_n^2 \leq \|\mathsf{f}_\lambda - f\|_n^2 + \sum_{j=1}^{M} r_{n,j}|\widehat{\lambda}_j - \lambda_j| + 2\sum_{j=1}^{M} r_{n,j}|\lambda_j| - 2\sum_{j=1}^{M} r_{n,j}|\widehat{\lambda}_j|.$$

Adding $\sum_{j=1}^{M} r_{n,j}|\widehat{\lambda}_j - \lambda_j|$ to both sides of this inequality and arguing as in (A.4), we get that, on the set $A$, for any $\lambda \in \mathbb{R}^M$,

$$\|\widetilde{f} - f\|_n^2 + \sum_{j=1}^{M} r_{n,j}|\widehat{\lambda}_j - \lambda_j| \leq \|\mathsf{f}_\lambda - f\|_n^2 + 4\sum_{j \in J(\lambda)} r_{n,j}|\widehat{\lambda}_j - \lambda_j|$$

$$\leq \|\mathsf{f}_\lambda - f\|_n^2 + 4\sqrt{M(\lambda)}\sqrt{\sum_{j \in J(\lambda)} r_{n,j}^2 |\widehat{\lambda}_j - \lambda_j|^2}.$$

Since $\sum\sum_{i \notin J(\lambda), j \notin J(\lambda)} \langle f_i, f_j\rangle_n (\widehat{\lambda}_i - \lambda_i)(\widehat{\lambda}_j - \lambda_j) \geq 0$ we have

$$\sum_{j \in J(\lambda)} r_{n,j}^2 |\widehat{\lambda}_j - \lambda_j|^2 = r_n^2 \|\widetilde{f} - \mathsf{f}_\lambda\|_n^2$$

$$- r_n^2 \sum_{i \notin J(\lambda), j \notin J(\lambda)} \sum \langle f_i, f_j\rangle_n (\widehat{\lambda}_i - \lambda_i)(\widehat{\lambda}_j - \lambda_j)$$

$$- 2r_n^2 \sum_{i \in J(\lambda), j > i} \sum \langle f_i, f_j\rangle_n (\widehat{\lambda}_i - \lambda_i)(\widehat{\lambda}_j - \lambda_j)$$

$$\leq r_n^2 \|\widetilde{f} - \mathsf{f}_\lambda\|_n^2 + 2r_n^2 \rho(\lambda) \left( \sum_{j=1}^M \|f_j\|_n |\widehat{\lambda}_j - \lambda_j| \right)^2.$$

Recalling that $r_{n,j} = r_n \|f_j\|_n$ and combining the last two displays, for $\lambda \in \mathbb{R}^M$ with $4\sqrt{2\rho(\lambda)M(\lambda)} \leq 1$, we obtain, on the set $A$,

$$\|\widetilde{f} - f\|_n^2 \leq \|\mathsf{f}_\lambda - f\|_n^2 + 4r_n\sqrt{M(\lambda)}(\|f - \mathsf{f}_\lambda\|_n + \|\widetilde{f} - f\|_n),$$

which is inequality (A.4) with $\kappa = 1/4$. The remainder of the proof now parallels that of Theorem 4.1. $\quad\square$

**A.5. Proof of Theorem 5.1.** We proceed similarly to [40]. The proof is based on the following easy corollary of the Fano lemma (which can be obtained, e.g., by combining Theorems 2.2 and 2.5 in [41]).

LEMMA A.1. *Let $w$ be a loss function, $A > 0$ be such that $w(A) > 0$, and let $\mathcal{C}$ be a set of functions on $\mathcal{X}$ of cardinality $N = \mathrm{card}(\mathcal{C}) \geq 2$ such that*

$$\|f - g\|_n^2 \geq 4s^2 > 0 \qquad \forall f, g \in \mathcal{C}, f \neq g,$$

*and the Kullback divergences $K(\mathbb{P}_f, \mathbb{P}_g)$ between the measures $\mathbb{P}_f$ and $\mathbb{P}_g$ satisfy*

$$K(\mathbb{P}_f, \mathbb{P}_g) \leq (1/16) \log N \qquad \forall f, g \in \mathcal{C}.$$

*Then for $\psi = s^2/A$ we have*

$$\inf_{T_n} \sup_{f \in \mathcal{C}} \mathbb{E}_f w[\psi^{-1} \|T_n - f\|_n^2] \geq c_1 w(A),$$

*where $\inf_{T_n}$ denotes the infimum over all estimators and $c_1 > 0$ is a constant.*

*The* (S) *aggregation case.* Pick $M$ disjoint subsets $S_1, \ldots, S_M$ of $\{X_1, \ldots, X_n\}$, each $S_j$ of cardinality $\log(M/D + 1)$ [w.l.o.g. we assume that $\log(M/D + 1)$ is an integer] and define the functions

$$f_j(x) = \gamma I_{\{x \in S_j\}}, \qquad j = 1, \ldots, M,$$

where $\gamma \leq L$ is a positive constant to be chosen. Consider the set of functions $\mathcal{V} = \{\mathsf{f}_\lambda : \lambda \in \bar{\Lambda}^{M,D}\}$ where $\bar{\Lambda}^{M,D}$ is the set of all $\lambda \in \mathbb{R}^M$ such that $D$ of the coordinates of $\lambda$ are equal to 1 and the remaining $M - D$ coordinates are zero. Clearly, $\mathcal{V} \subset \mathcal{F}_0$. Thus, it suffices to prove the (S) lower bound of the theorem where the supremum over $f \in \mathcal{F}_0$ is replaced by that over $f \in \mathcal{V}$. Since $\bar{\Lambda}^{M,D} \subset \Lambda^{M,D}$, for $f \in \mathcal{V}$ we have $\min_{\lambda \in \Lambda^{M,D}} \|\mathsf{f}_\lambda - f\|_n^2 = 0$. Therefore, to finish the proof for the (S) case, it suffices to bound from below the quantity $\inf_{T_n} \sup_{f \in \mathcal{V}} \mathbb{E}_f w(\psi_{n,M}^{-1} \|T_n - f\|_n^2)$ where $\psi_{n,M} = D \log(M/D +$ $1)/n$. This will be done by applying Lemma A.1. In fact, note that for every two functions $f_\lambda$ and $f_{\bar\lambda}$ in $\mathcal{V}$ we have

$$(A.6) \qquad \|f_\lambda - f_{\bar\lambda}\|_n^2 = \frac{\gamma^2 \log(M/D+1)}{n} \rho(\lambda, \bar\lambda) \le \frac{\gamma^2 D \log(M/D+1)}{n},$$

where $\rho(\lambda, \bar\lambda) \overset{\text{def}}{=} \sum_{j=1}^M I_{\{\lambda_j \neq \bar\lambda_j\}}$ is the Hamming distance between $\lambda = (\lambda_1, \ldots, \lambda_M) \in \bar\Lambda^{M,D}$ and $\bar\lambda = (\bar\lambda_1, \ldots, \bar\lambda_M) \in \bar\Lambda^{M,D}$. Lemma 4 in [10] (see also [22]) asserts that if $M \ge 6D$ there exists a subset $\Lambda' \subset \Lambda^{M,D}$ such that, for some constant $\tilde c > 0$ independent of $M$ and $D$,

$$(A.7) \qquad \log \operatorname{card}(\Lambda') \ge \tilde c D \log\left(\frac{M}{D} + 1\right)$$

and

$$(A.8) \qquad \rho(\lambda, \bar\lambda) \ge \tilde c D \qquad \forall\, \lambda, \bar\lambda \in \Lambda', \lambda \neq \bar\lambda.$$

Consider a set of functions $\mathcal{C} = \{f_\lambda : \lambda \in \Lambda'\} \subset \mathcal{V}$. From (A.6) and (A.8), for any two functions $f_\lambda$ and $f_{\bar\lambda}$ in $\mathcal{C}$ we have

$$(A.9) \qquad \|f_\lambda - f_{\bar\lambda}\|_n^2 \ge \frac{\tilde c \gamma^2 D \log(M/D+1)}{n} \overset{\text{def}}{=} 4s^2.$$

Since the $W_j$'s are $N(0, \sigma^2)$ random variables, the Kullback divergence $K(\mathbb{P}_{f_\lambda}, \mathbb{P}_{f_{\bar\lambda}})$ between $\mathbb{P}_{f_\lambda}$ and $\mathbb{P}_{f_{\bar\lambda}}$ satisfies

$$(A.10) \qquad K(\mathbb{P}_{f_\lambda}, \mathbb{P}_{f_{\bar\lambda}}) = \frac{n}{2\sigma^2} \|f_\lambda - f_{\bar\lambda}\|_n^2, \qquad j = 1, \ldots, M.$$

In view of (A.6) and (A.10), one can choose $\gamma$ small enough to have

$$K(\mathbb{P}_{f_\lambda}, \mathbb{P}_{f_{\bar\lambda}}) \le \frac{1}{16\tilde c} D \log\left(\frac{M}{D} + 1\right) \le \frac{1}{16} \log \operatorname{card}(\Lambda') = \frac{1}{16} \log \operatorname{card}(\mathcal{C})$$

for all $\lambda, \bar\lambda \in \Lambda'$. Now, to get the lower bound for the (S) case, it remains to use this inequality together with (A.9) and to apply Lemma A.1. Thus, the (S) lower bound is proved under the assumption that $M \ge 6D$, which is needed to assure (A.7) and (A.8).

In the remaining case where $D \le M < 6D$ we use another construction. Note that it is enough to prove the result for $\psi_{n,M} = D/n$. We consider separately the cases $D \ge 8$ and $2 \le D < 8$. If $D \ge 8$ we consider the functions $f_j(x) = \gamma I_{\{x = X_j\}}$, $j = 1, \ldots, M$, and a finite set of their linear combinations,

$$(A.11) \qquad \mathcal{U} = \left\{ g = \sum_{j=1}^D \omega_j f_j : \omega \in \Omega \right\},$$

where $\Omega$ is the set of all vectors $\omega \in \mathbb{R}^M$ with binary coordinates $\omega_j \in \{0, 1\}$. Since the supports of the $f_j$'s are disjoint, the functions $g \in \mathcal{U}$ are uniformly bounded by $\gamma$, thus $\mathcal{U} \subset \mathcal{F}_0$. Also, $\mathcal{U} \subset \{f_\lambda : \lambda \in \Lambda^{M,D}\}$ since at most the first $D$ functions $f_j$ are included in the linear combination. Clearly, $\min_{\lambda \in \Lambda^{M,D}} \|f_\lambda - f\|_n^2 = 0$ for any $f \in \mathcal{U}$. Therefore it remains to bound from below the quantity $\inf_{T_n} \sup_{f \in \mathcal{U}} \mathbb{E}_f w(\psi_{n,M}^{-1} \|T_n - f\|_n^2)$, where $\psi_{n,M} = D/n$. To this end, we apply again Lemma A.1.

Note that for any $g_1 = \sum_{j=1}^D \omega_j f_j \in \mathcal{U}$ and $g_2 = \sum_{j=1}^D \omega_j' f_j \in \mathcal{U}$ we have

$$(A.12) \qquad \|g_1 - g_2\|_n^2 = \frac{\gamma^2}{n} \sum_{j=1}^D (\omega_j - \omega_j')^2 \leq \gamma^2 D/n.$$

Since $D \geq 8$ it follows from the Varshamov–Gilbert bound (see [22] or [41], Chapter 2) that there exists a subset $\mathcal{C}'$ of $\mathcal{U}$ such that $\mathrm{card}(\mathcal{U}_0) \geq 2^{D/8}$ and

$$(A.13) \qquad \|g_1 - g_2\|_n^2 \geq C_1 \gamma^2 D/n$$

for any $g_1, g_2 \in \mathcal{C}'$. Using (A.10) and (A.12) we get, for any $g_1, g_2 \in \mathcal{C}'$,

$$K(\mathbb{P}_{g_1}, \mathbb{P}_{g_2}) \leq C_2 \gamma^2 D \leq C_3 \gamma^2 \log(\mathrm{card}(\mathcal{C}')),$$

and choosing $\gamma$ small enough, we can finish the proof by applying Lemma A.1 where we take $\mathcal{C} = \mathcal{C}'$ and act in the same way as above for $M \geq 6D$.

Finally, if $D \leq M < 6D$ and $2 \leq D < 8$, we have $\psi_{n,M} < 8/n$, and the proof is easily obtained by choosing $f_1 \equiv 0$ and $f_2 \equiv \gamma n^{-1/2}$ and applying Lemma A.1 to the set $\mathcal{C} = \{f_1, f_2\}$.

*The* (MS) *aggregation case.* We use the proof for (S) aggregation with $D = 1$. Note that $\bar{\Lambda}^{M,1}$ is the set of all the vertices of $\Lambda^M$, except the vertex $(0, \dots, 0)$. Thus, the proof for the (S) case with $M \geq 6D$ and $D = 1$ gives us the required lower bound for the (MS) case, with the optimal rate $\psi_{n,M} = (\log M)/n$. It remains to treat (MS) aggregation for $M < 6$. Then we have $\psi_{n,M} \leq (\log 7)/n$, and we apply Lemma A.1 to the set $\mathcal{C} = \{f_{\lambda'}, f_{\lambda''}\}$ where $\lambda' = (1, 0, \dots, 0) \in \Lambda^M$, $\lambda'' = (0, \dots, 0, 1) \in \Lambda^M$ and $f_\lambda$ is defined in the proof for the (S) case. Clearly, $\|f_{\lambda'} - f_{\lambda''}\|_n^2 = 2\gamma^2 \log(M+1)/n \geq 2\gamma^2 (\log 3)/n$, and the result easily follows from (A.10) and Lemma A.1.

*The* (C) *aggregation case.* Consider the orthonormal trigonometric basis in $L_2[0,1]$ defined by $\phi_1(x) \equiv 1$, $\phi_{2k}(x) = \sqrt{2} \cos(2\pi k x)$, $\phi_{2k+1}(x) = \sqrt{2} \sin(2\pi \times kx)$, $k = 1, 2, \dots$, for $x \in [0,1]$. Set

$$(A.14) \qquad f_j(x) = \gamma \sum_{k=1}^n \phi_j(k/n) I_{\{x = X_k\}}, \qquad j = 1, \dots, M,$$

where $\gamma \leq L/\sqrt{2}$ is a positive constant to be chosen. The system of functions $\{\phi_j\}_{j=1,\dots,M}$ is orthonormal w.r.t. the discrete measure that assigns mass $1/n$ to each of the points $k/n, k = 1, \dots, n$:

$$\frac{1}{n} \sum_{k=1}^n \phi_j(k/n)\phi_l(k/n) = \delta_{jl}, \qquad j, l = 1, \dots, n,$$

where $\delta_{jl}$ is the Kronecker delta (see, e.g., [41], Lemma 1.9). Hence

(A.15) $$\langle f_j, f_l \rangle_n = \gamma^2 \delta_{jl}, \qquad j, l = 1, \ldots, M,$$

where $\langle \cdot, \cdot \rangle_n$ stands for the scalar product induced by $\| \cdot \|_n$.

Assume first that $M > \sqrt{n}$ (i.e., we are in the "sparse" case). Define an integer

$$m = \left\lceil c_2 \left[ n \Big/ \log\left( \frac{M}{\sqrt{n}} + 1 \right) \right]^{1/2} \right\rceil$$

for a constant $c_2 > 0$ chosen in such a way that $M \geq 6m$. Consider the finite set $\mathcal{C} \subset \Lambda^M$ composed of such convex combinations of $f_1, \ldots, f_M$ that $m$ of the coefficients $\lambda_j$ are equal to $1/m$ and the remaining $M - m$ coefficients are zero. In view of (A.15), for every pair of functions $g_1, g_2 \in \mathcal{C}$ we have

(A.16) $$\|g_1 - g_2\|_n^2 \leq 2\gamma^2/m.$$

To finish the proof for $M > \sqrt{n}$ it suffices now to apply line-by-line the argument after the formula (10) in [40] replacing there $\| \cdot \|$ by $\| \cdot \|_n$. Similarly, the proof for $M \leq \sqrt{n}$ is analogous to that given in [40], with the only difference that the functions $f_j$ should be chosen as in (A.14) and $\| \cdot \|$ should be replaced by $\| \cdot \|_n$.

## APPENDIX B: TECHNICAL LEMMAS

LEMMA B.1. *Let $f, f_1, \ldots, f_M \in \mathcal{F}_0$ and $1 \leq m \leq M$. Let $\mathcal{C}$ be the finite set of functions defined in the proof of (3.5). Then (A.1) holds:*

(B.1) $$\min_{g \in \mathcal{C}} \|g - f\|_n^2 \leq \min_{\lambda \in \Lambda^M} \|\mathsf{f}_\lambda - f\|_n^2 + \frac{L^2}{m}.$$

PROOF. Let $f^*$ be a minimizer of $\|\mathsf{f}_\lambda - f\|_n^2$ over $\lambda \in \Lambda^M$. Clearly, $f^*$ is of the form

$$f^* = \sum_{j=1}^{M} p_j f_j \qquad \text{with } p_j \geq 0 \text{ and } \sum_{j=1}^{M} p_j \leq 1.$$

Define a probability distribution on $j = 0, 1, \ldots, M$ by

$$\pi_j = \begin{cases} p_j, & \text{if } j \neq 0, \\ 1 - \sum_{j=1}^{M} p_j, & \text{if } j = 0. \end{cases}$$

Consider $m$ i.i.d. random integers $j_1, \ldots, j_m$ where each $j_k$ is distributed according to $\{\pi_j\}$ on $\{0, 1, \ldots, M\}$. Introduce the random function

$$\bar{f}_m = \frac{1}{m} \sum_{k=1}^{m} g_{j_k}, \qquad \text{where } g_j = \begin{cases} f_j, & \text{if } j \neq 0, \\ 0, & \text{if } j = 0. \end{cases}$$

For every $x \in \mathcal{X}$ the random variables $g_{j_1}(x), \ldots, g_{j_m}(x)$ are i.i.d. with $\mathbb{E}(g_{j_k}(x)) = f^*(x)$. Thus,

$$\mathbb{E}(\bar{f}_m(x) - f^*(x))^2 = \mathbb{E}\left(\left[\frac{1}{m}\sum_{k=1}^m \{g_{j_k}(x) - \mathbb{E}(g_{j_k}(x))\}\right]^2\right)$$

$$\leq \frac{1}{m}\mathbb{E}(g_{j_1}^2(x)) \leq \frac{L^2}{m}.$$

Hence for every $x \in \mathcal{X}$ and every $f \in \mathcal{F}_0$ we get

$$(B.2) \qquad \begin{aligned} \mathbb{E}(\bar{f}_m(x) - f(x))^2 &= \mathbb{E}(\bar{f}_m(x) - f^*(x))^2 + (f^*(x) - f(x))^2 \\ &\leq \frac{L^2}{m} + (f^*(x) - f(x))^2. \end{aligned}$$

Integrating (B.2) over the empirical probability measure that puts mass $1/n$ at each $X_i$ and recalling the definition of $f^*$, we obtain

$$(B.3) \qquad \mathbb{E}\|\bar{f}_m - f\|_n^2 \leq \min_{\lambda \in \Lambda^M} \|f_\lambda - f\|_n^2 + \frac{L^2}{m}.$$

Finally, note that the random function $\bar{f}_m$ takes its values in $\mathcal{C}$, which implies that

$$\mathbb{E}\|\bar{f}_m - f\|_n^2 \geq \min_{g \in \mathcal{C}} \|g - f\|_n^2.$$

This and (B.3) prove (B.1).  □

LEMMA B.2.  *Let $Z_d$ denote a random variable having the $\chi^2$ distribution with $d$ degrees of freedom. Then for all $x > 0$,*

$$(B.4) \qquad \mathbb{P}\{Z_d - d \geq x\sqrt{2d}\} \leq \exp\left(-\frac{x^2}{2(1 + x\sqrt{2/d})}\right).$$

PROOF.  See [16], page 857.  □

**Acknowledgments.**  We would like to thank the Associate Editor, the referees and Lucien Birgé for insightful comments.

F. Bunea
M. H. Wegkamp
Department of Statistics
Florida State University
Tallahassee, Florida
USA
E-mail: bunea@stat.fsu.edu
          wegkamp@stat.fsu.edu

A. B. Tsybakov
Laboratoire de Probabilités
  et Modèles Aléatoires
Université Paris VI
France
E-mail: tsybakov@ccr.jussieu.fr